\documentclass[12pt]{article}
\usepackage{amsmath}
\usepackage{hyperref}
\usepackage{amsfonts}
\usepackage{amssymb}
\usepackage{graphics}
\usepackage{mathdots}
\usepackage{wrapfig}
\usepackage{color}
\usepackage[pdftex]{graphicx}

\def\bee{\begin{equation}}
\def\eee{\end{equation}}

\def\DHLhksqrt#1#2{\setbox0=\hbox{$#1\sqrt{#2\,}$}\dimen0=\ht0
\advance\dimen0-0.2\ht0
\setbox2=\hbox{\vrule height\ht0 depth -\dimen0}%
{\box0\lower0.4pt\box2}}
\allowdisplaybreaks

\begin{document}

\bigskip\bigskip\bigskip

\bigskip
\centerline{    }
\centerline{\Large\bf  On the moments of the gaps between consecutive primes}
\bigskip\bigskip
\centerline{\large\sl Marek Wolf}
\begin{center}
Cardinal  Stefan  Wyszynski  University, Faculty  of Mathematics and Natural Sciences. College of Sciences,\\
ul. W{\'o}ycickiego 1/3,   PL-01-938 Warsaw,   Poland, e-mail:  m.wolf@uksw.edu.pl
\end{center}

\bigskip\bigskip

\begin{abstract}
We derive heuristically   formula for the $k$--moments  $M_k(x)$  of the gaps between consecutive primes$<x $  represented directly
by  $x$$\pi(x)$ --- the number of primes up to:  $M_k(x)= \Gamma(k+1)x^k/\pi^{k-1}(x)+\mathcal{O}(x)$, We  illustrate obtained  results  by  computer data.
\end{abstract}
\bigskip\bigskip
\bigskip\bigskip



\bibliographystyle{abbrv}

Let $p_n$ denotes the $n$-th prime number and $d_n=p_{n+1}-p_n$ denotes the $n$-th gap between consecutive primes.
Let us introduce  moments of arbitrary order $k$  of gaps between consecutive primes:
\bee
M_k(x) \equiv \sum_{p_{n+1}\leq x}(p_{n+1}-p_n)^k.
\eee
The   symbol $f(x)\sim  g(x)$ means here that $\lim_{x\rightarrow \infty} f(x)/g(x)=1$. Presumably for the first time the
second moment of gaps $M_2(x)\sum_{p_n<x} (p_{n+1}-p_n)^2$  was considered in 1937
by H. Cramer \cite{Cramer}. Assuming the validity of the Riemann Hypothesis he obtained:
\bee
\sum_{p_n<x} (p_{n+1}-p_n)^2 =  \mathcal{O}(x\log^{3+\epsilon} x)
\eee
for every $\epsilon >0$. In 1943    A. Selberg in \cite{Selberg1943}, also  assuming the Riemann Hypothesis,  has proved:
\bee
\sum_{p_n<x} \frac{(p_{n+1}-p_n)^2}{p_n} =  \mathcal{O}(\log^{3} x).
\eee
In \cite{Heath-Brown}  D.R.  Heath-Brown  conjectured that
\bee
\sum_{p_{n+1}<x}(p_{n+1}-p_n)^2 \sim 2x \log(x).
\label{eq-Heath-Brown-2}
\eee
For the history of the problem and review of results see \cite{Heath-Brown-1990};  see also problem A8 in \cite{Guy}.
In \cite[p.2056]{Oliveira-2013}  the Heath-Brown--Oliveira   conjecture was formulated:
\bee
M_k(x)  \sim k! x\log^{k-1} (x).
\label{eq-Heath-Brown}
\eee
In \cite{Oliveira-2013} authors made a remark after   equation \eqref{eq-Heath-Brown} that $k\geq 1$,  but even for $k=0$ it produces correct
answer as $M_0(x)$  is by 1 less then  the number of primes up to $x$:  $M_0(x)=\pi(x) -1 $ (here, as  usual,
$\pi(x)=\sum_{n} \Theta(x-p_n)$ and  $\Theta$ is a unit step function: $\Theta(x)=1$  for $x>0$ and $\Theta(x)=0$ for $x\leq0$).
By the Prime Number Theorem (PNT) the  number of prime numbers  below $x$ is very well approximated by the logarithmic integral
\[
\pi(x)\sim  {\rm Li}(x) \equiv \int_2^x \frac{du}{\ln(u)}.
\]
Integration by parts gives the asymptotic expansion which should be cut at the term $n_0=\lfloor \ln(x)\rfloor$:
\bee
{\rm Li}(x) = \frac{x}{\ln(x)}  + \frac{x}{\ln^2(x)} + \frac{2!x}{\ln^3(x)} + \frac{3!x}{\ln^4(x)} +\cdots .
\label{PNT}
\eee

\bigskip

Let $\tau_d(x)$ denote the number of pairs of  consecutive   primes smaller than a given bound $x$ and separated by $d$:
\bee
\tau_d(x)= \sharp\{ p_n,  p_{n+1} < x,~~{\rm  with}~ p_{n+1}-p_n=d\}.
\label{definition}
\eee

In \cite{Wolf-heuristics} (see also  \cite{Wolf-PRE}) we proposed the following formula expressing function $\tau_d(x)$
directly by $\pi(x)$:
\bee
\tau_d(x) \sim  C_2 \prod_{p \mid d, p > 2} \frac{p - 1}{p - 2}~~ \frac{\pi^2(x)}{x}\Bigg(1-\frac{2\pi(x)}{x}\Bigg)^{{\frac{ d}{2}-1}} ~~~{\rm for} ~ d\geq 6.
\label{main}
\eee
Here
\[C_2 \equiv 2\prod_{p > 2} \biggl( 1 - \frac{1}{(p - 1)^2}\biggr) =1.320323631693739\ldots
\]
is  called the ``twins constant''. The pairs of primes separated by $d=2$ (``twins'') and $d=4$  (``cousins'')  are special as they
always have to be consecutive primes (with the exception of the pair  (3,7) containing 5 in the middle)).   For $d=4$ we
adapt the expression obtained from \eqref{main} for $d=2$, which for $\pi(x)\sim x/\log(x)$  goes into the
the  conjecture B of G. H. Hardy and J.E. Littlewood  \cite[eqs. (5.311) and (5.312)]{Hardy_and_Littlewood}:
\bee
\tau_2(x)\big(\approx \tau_4(x)\big)\sim  C_2 \frac{\pi^2(x)}{x} \approx C_2 \frac{x}{\ln^2(x)}.
\eee

We will assume that for sufficiently regular functions  $f(n)$ the  following  formula  holds:
\bee
\sum_{k=1}^\infty \prod_{p\mid k, p>2}{p-1\over p-2} ~f(k) = {1 \over \prod_{p
> 2}( 1 - {1 \over (p - 1)^2})} \sum_{k=1}^\infty f(k)
\label{rownosc1}
\eee
In other words we will  replace the product over $p|d$ in \eqref{main} by its mean value as  E. Bombieri and H. Davenport
\cite{Bombieri} have proved that the number  $1/\prod_{p > 2}( 1 - {1 \over (p - 1)^2})C_2/2$ is the arithmetical average
of the product $\prod_{p\mid k} \frac{p-1}{p-2}$:
\bee
\sum_{k=1}^n \prod_{p\mid k,p>2}{p-1\over p-2} =
{n \over \prod_{p > 2}( 1 - {1 \over (p - 1)^2})} + \mathcal{O}(\log^2(n)).
\eee
Later H.L. Montgomery \cite[eq.(17.11)]{Montgomery} has improved the error term to $\mathcal{O}(\log(n))$.

\vskip 0.4cm
\begin{center}
{\sf TABLE {\bf I}} ~The ratios of the sums of squares of gaps between consecutive primes
for $x=2^{24}\ldots, x=4\ldots  \times 10^{18}$ and closed formulas for $M_2(x)$  given by
eq.(\ref{eq-Heath-Brown}),  eq.(\ref{prediction1}) and eq.(\ref{eq-moment-2}) respectively  presented  up to 4 figures.
\begin{tabular}{|c|c|c|c|} \hline
$x$ & $ M_2(x) /\widetilde{M}_2^{(1)}(x)$ & $  M_2(x)/\widetilde{M}_2^{(2)}(x) $ &   $  M_2(x) /\widetilde{M}_2^{(3)}(x)$\\ \hline
$2^{24}=1.6777\!  \ldots  \times \! 10^{7 }$ &      0.7971&    0.9104 &    0.8519 \\ \hline
$2^{26}=6.7109\!  \ldots  \times \! 10^{7 }$ &      0.8102&    0.9151 &    0.8611 \\ \hline
$2^{28}=2.6844\!  \ldots  \times \! 10^{8 }$ &      0.8221&    0.9198 &    0.8696 \\ \hline
$2^{30}=1.0737\!  \ldots  \times \! 10^{9 }$ &      0.8323&    0.9237 &    0.8769 \\ \hline
$2^{32}=4.2950\!  \ldots  \times \! 10^{9 }$ &      0.8414&    0.9272 &    0.8833 \\ \hline
$2^{34}=1.7180\!  \ldots  \times \! 10^{10 }$ &      0.8495&    0.9303 &    0.8890 \\ \hline
$2^{36}=6.8719\!  \ldots  \times \! 10^{10 }$ &      0.8567&    0.9332 &    0.8942 \\ \hline
$2^{38}=2.7488\!  \ldots  \times \! 10^{11 }$ &      0.8632&    0.9358 &    0.8988 \\ \hline
$2^{40}=1.0995\!  \ldots  \times \! 10^{12 }$ &      0.8692&    0.9382 &    0.9031 \\ \hline
$2^{42}=4.3980\!  \ldots  \times \! 10^{12 }$ &      0.8746&    0.9404 &    0.9069 \\ \hline
$2^{44}=1.7592\!  \ldots  \times \! 10^{13 }$ &      0.8796&    0.9425 &    0.9105 \\ \hline
$2^{46}=7.0369\!  \ldots  \times \! 10^{13 }$ &      0.8841&    0.9444 &    0.9138 \\ \hline
$2^{48}=2.8147\!  \ldots  \times \! 10^{14 }$ &      0.8883&    0.9462 &    0.9168 \\ \hline
$1.61    \times 10^{18}$ &      0.9087&    0.9549 &    0.9315 \\ \hline
$4    \times 10^{18}$ &      0.9104&    0.9556 &    0.9327 \\ \hline
\end{tabular} \\
\end{center}
\vskip 0.4cm

We will use the notation $\widetilde{M}_k^{(i)}(x)$  for the  $i$-th analytical formula for $M_k(x)$.  The superscript $i=1$ will
refer to the  conjecture \eqref{eq-Heath-Brown}:   $\widetilde{M}_k^{(1)}(x)=k!x\log^{k-1}(x)$  and expressions for $i=2$ and $i=3$
we will  derive below.   For second moments  using
the differentiated geometrical series  we obtain (we have extended the summation over $d=2n$
 up to infinity and used \eqref{rownosc1},  then the dependence on $c_2$ drops out)
\[
M_2(x)  =\sum_{p_n<x} (p_n-p_{n-1})^2 = \sum_{d=2, 4,6,\ldots} d^2 \tau_d(x)\approx  \frac{8\pi^2(x)}{(x-2\pi(x))}\sum_{n=1}^\infty n^2\Big(1-\frac{2\pi(x)}{x}\Big)^n
\]
\bee
=\frac{2x^2}{\pi(x)}\Big(1-\frac{\pi(x)}{x}\Big) \equiv \widetilde{M}_2^{(2)}(x).
\label{prediction1}
\eee
For large $x$ skipping in the big bracket   above term   $\pi(x)/x\sim 1/\log(x)$
we obtain
\bee
M_2(x)  \sim \frac{2 x^2}{\pi(x)}\equiv \widetilde{M}_2^{(3)}(x)
\label{eq-moment-2}
\eee
what for $\pi(x) \sim x/\log(x)$ gives exactly (\ref{eq-Heath-Brown-2}).

In the similar manner for third moment we obtain using \eqref{main} the expression:
\bee
M_3(x)\sim \widetilde{M}_3^{(2)}(x) \equiv  \frac{6x^3}{\pi^2(x)}\Big(1-2\frac{\pi(x)}{x}+\frac{2}{3}\frac{\pi^2(x)}{x^2}\Big).    
\label{eq-moment-3}
\eee
Putting  here $\pi(x) \sim x/\log(x)$  in the limit of large $x$  we obtain $M_3(x)\sim 6 x \log^2(x)$,  i.e. \eqref{eq-Heath-Brown-2}
for  $k=3$.

For fourth moment similarly we obtain:
\bee
M_4(x)\approx \widetilde{M}_4^{(2)}(x)  \equiv 24\frac{x^4}{\pi^3(x)}\Big(1-3\frac{\pi(x)}{x}+\frac{7}{3}\frac{\pi^2(x)}{x^2}-\frac{1}{3}\frac{\pi^3(x)}{x^3}\Big).
\label{eq-moment-4}
\eee
and for large $x$  it goes to $4! x\log^3(x)$.

\vskip 0.4cm
\begin{center}
{\sf TABLE {\bf II}}   ~The ratios of the sums of cubes  of gaps between consecutive primes
for $x=2^{24}\ldots, x=4\ldots  \times 10^{18}$ and closed formulas for $M_3(x)$  given by eq.(\ref{eq-Heath-Brown}) for $k=3$,
 eq.(\ref{eq-moment-3}) and  eq.(\ref{prediction2}) for $k=3$  presented  up to 4 figures.
\begin{tabular}{|c|c|c|c|} \hline
$x$ & $ M_3(x) /\widetilde{M}_3^{(1)}(x)$ & $  M_3(x)/\widetilde{M}_3^{(2)}(x) $ &   $  M_3(x) /\widetilde{M}_3^{(3)}(x)$\\ \hline
$2^{24}=1.6777\! \! \ldots  \times \! 10^{  7}$ &   0.6104&    0.7975 &    0.6972 \\ \hline
$2^{26}=6.7109\! \! \ldots  \times \! 10^{  7}$ &   0.6331&    0.8087 &    0.7152 \\ \hline
$2^{28}=2.6844\! \! \ldots  \times \! 10^{  8}$ &   0.6540&    0.8195 &    0.7318 \\ \hline
$2^{30}=1.0737\! \! \ldots  \times \! 10^{  9}$ &   0.6722&    0.8287 &    0.7461 \\ \hline
$2^{32}=4.2950\! \! \ldots  \times \! 10^{  9}$ &   0.6885&    0.8367 &    0.7588 \\ \hline
$2^{34}=1.7180\! \! \ldots  \times \! 10^{ 10}$ &   0.7030&    0.8438 &    0.7700 \\ \hline
$2^{36}=6.8719\! \! \ldots  \times \! 10^{ 10}$ &   0.7162&    0.8504 &    0.7803 \\ \hline
$2^{38}=2.7488\! \! \ldots  \times \! 10^{ 11}$ &   0.7283&    0.8564 &    0.7896 \\ \hline
$2^{40}=1.0995\! \! \ldots  \times \! 10^{ 12}$ &   0.7393&    0.8619 &    0.7981 \\ \hline
$2^{42}=4.3980\! \! \ldots  \times \! 10^{ 12}$ &   0.7495&    0.8670 &    0.8059 \\ \hline
$2^{44}=1.7592\! \! \ldots  \times \! 10^{ 13}$ &   0.7588&    0.8716 &    0.8131 \\ \hline
$2^{46}=7.0369\! \! \ldots  \times \! 10^{ 13}$ &   0.7674&    0.8759 &    0.8198 \\ \hline
$2^{48}=2.8147\! \! \ldots  \times \! 10^{ 14}$ &   0.7754&    0.8800 &    0.8259 \\ \hline
$1.61   \times 10^{18}$ &       0.8147&    0.8997 &    0.8561 \\ \hline
$4  \times 10^{18}$ &       0.8180&    0.9014 &    0.8586 \\ \hline
\end{tabular}
\end{center}
\vskip 0.4cm

\vskip 0.4cm
\begin{center}
{\sf TABLE {\bf III}}   ~The ratios of the sums of fourth powers  of gaps between consecutive primes
for $x=2^{24}\ldots, x=4\ldots  \times 10^{18}$ and closed formulas for $M_4(x)$  given by eq.(\ref{eq-Heath-Brown}) for $k=4$,
 eq.(\ref{eq-moment-4}) and  eq.(\ref{prediction2}) for $k=4$ presented up to 4 figures.
\begin{tabular}{|c|c|c|c|} \hline
$x$ & $ M_4(x) /\widetilde{M}_4^{(1)}(x)$ & $  M_4(x)/\widetilde{M}_4^{(2)}(x) $ &   $  M_4(x) /\widetilde{M}_4^{(3)}(x)$\\ \hline
$2^{24}=1.6777\!  \ldots  \times \! 10^{  7}$ &      0.4586 &    0.6854 &    0.5598 \\ \hline
$2^{26}=6.7109\!  \ldots  \times \! 10^{  7}$ &      0.4862 &    0.7024 &    0.5838 \\ \hline
$2^{28}=2.6844\!  \ldots  \times \! 10^{  8}$ &      0.5123 &    0.7190 &    0.6063 \\ \hline
$2^{30}=1.0737\!  \ldots  \times \! 10^{  9}$ &      0.5354 &    0.7332 &    0.6261 \\ \hline
$2^{32}=4.2950\!  \ldots  \times \! 10^{  9}$ &      0.5560 &    0.7453 &    0.6433 \\ \hline
$2^{34}=1.7180\!  \ldots  \times \! 10^{ 10}$ &      0.5746 &    0.7560 &    0.6587 \\ \hline
$2^{36}=6.8719\!  \ldots  \times \! 10^{ 10}$ &      0.5919 &    0.7661 &    0.6731 \\ \hline
$2^{38}=2.7488\!  \ldots  \times \! 10^{ 11}$ &      0.6078 &    0.7753 &    0.6861 \\ \hline
$2^{40}=1.0995\!  \ldots  \times \! 10^{ 12}$ &      0.6225 &    0.7837 &    0.6982 \\ \hline
$2^{42}=4.3980\!  \ldots  \times \! 10^{ 12}$ &      0.6360 &    0.7915 &    0.7093 \\ \hline
$2^{44}=1.7592\!  \ldots  \times \! 10^{ 13}$ &      0.6486 &    0.7987 &    0.7195 \\ \hline
$2^{46}=7.0369\!  \ldots  \times \! 10^{ 13}$ &      0.6603 &    0.8054 &    0.7290 \\ \hline
$2^{48}=2.8147\!  \ldots  \times \! 10^{ 14}$ &      0.6712 &    0.8116 &    0.7379 \\ \hline
$1.61  \times 10^{18}$ &       0.7256&    0.8422 &    0.7816 \\ \hline
$4   \times 10^{18}$ &       0.7303&    0.8448 &    0.7853 \\ \hline
\end{tabular}
\end{center}
\vskip 0.4cm

We stop with these particular  moments and we will derive the formula for  moments  of  general order $k$.
From the  formula (\ref{main}) we obtain  :\\
\bee
M_k(x)=\sum_{p_n<x} (p_n-p_{n-1})^k = \sum_{d=2, 4,6,\ldots} d^k \tau_d(x) \sim 2\frac{\pi^2(x)}{x-2\pi(x)}\sum_{n=1}^\infty (2n)^k \Big(1-\frac{2\pi(x)}{x}\Big)^n
\eee
To proceed further we need  formula for the $k$-times differentiated  geometrical series:
\bee
\sum_{n=1}^\infty  n^kq^n = \Big(q\frac{d}{dq}\Big)^k\frac{1}{1-q} = \frac {1}{(1-q)^{k+1}}\sum_{i=0}^{k-1} \genfrac{\langle}{\rangle}{0pt} {}{ k}{ i}  q^{k-i},
\eee
where  $|q|<1$ and $\genfrac{\langle}{\rangle}{0pt} {}{ n}{ i}$   are  Eulerian numbers (should not be confused with Euler
numbers $E_n$),  see \cite[p.54]{Petersen-book} and eq. (7) in entry {\it   Eulerian numbers} in \cite{Weisstein}.
In our case $q=1-2\pi(x)/x$  and for large  $x$  we have $q\to 1$   hence in nominator we obtain $k!$ because the
Eulerian numbers  satisfy the  identity
\bee
\sum_{n=0}^k\genfrac{\langle}{\rangle}{0pt} {}{ k}{ n} =k!,
\eee
see  \cite[eq.(1.8)]{Carlitz-1972} and entry {\it   Eulerian numbers} in  \cite{Weisstein}. The denomiator
is  $(2\pi(x)/x)^{k+1}$ and the power $2^{k+1}$ cancels out.  Finally we obtain
\bee
M_k(x)\equiv \sum_{p_n<x} (p_n-p_{n-1})^k  \sim   k! \frac{x^k}{\pi^{k-1}(x)} \equiv  \widetilde{M}_k^{(3)}(x)
\label{prediction2}
\eee
and for $\pi(x)\sim x/\log(x)$ it goes into  \eqref{eq-Heath-Brown}. For $k=1$ from above equation we obtain $M_1(x)=x$
and for $k=0$ we obtain $M_0(x)=\pi(x)$ as it should be.

During over a seven months long run of the computer  program we have collected  the values of $\tau_d(x)$ up to
$x=2^{48}\approx 2.8147\times 10^{14}$.  The data representing the function $\tau_d(x)$ were stored at values of $x$
forming the geometrical  progression with the ratio 2, i.e. at $x=2^{15}, 2^{16}, \ldots, 2^{47},
2^{48}$. Such a choice of the intermediate thresholds as powers of 2  was determined by the employed computer program
in which  the primes were coded  as bits.   The data is available for downloading from  \url{http://pracownicy.uksw.edu.pl/mwolf/gaps.zip}.
At the Tom{\'a}s   Oliveira e Silva web site \url{http://sweet.ua.pt/tos/gaps.html}  we have found values of $\tau_d(x)$ for
$x=1.61\times 10^{18}$  and  $x=4\times 10^{18}$.   In the tables  I, II and III we  present comparison of the actual  values of
$M_k(x)$  calculated  from these computer data  $k=2, 3, 4$  and the prediction given by  formulas for  $\widetilde{M}_k^{(i)}(x)$
for $i=1,2,3$  and the set of values of $x$.  As the rule the best approximations are   given by  \eqref{eq-moment-2},
\eqref{eq-moment-3}  and \eqref{eq-moment-4},  next by  \eqref{prediction2} and the least accurate  are values predicted by
\eqref{eq-Heath-Brown}.

\begin{figure*}[t]
\centering
\includegraphics[width=0.95\textwidth, angle=0]{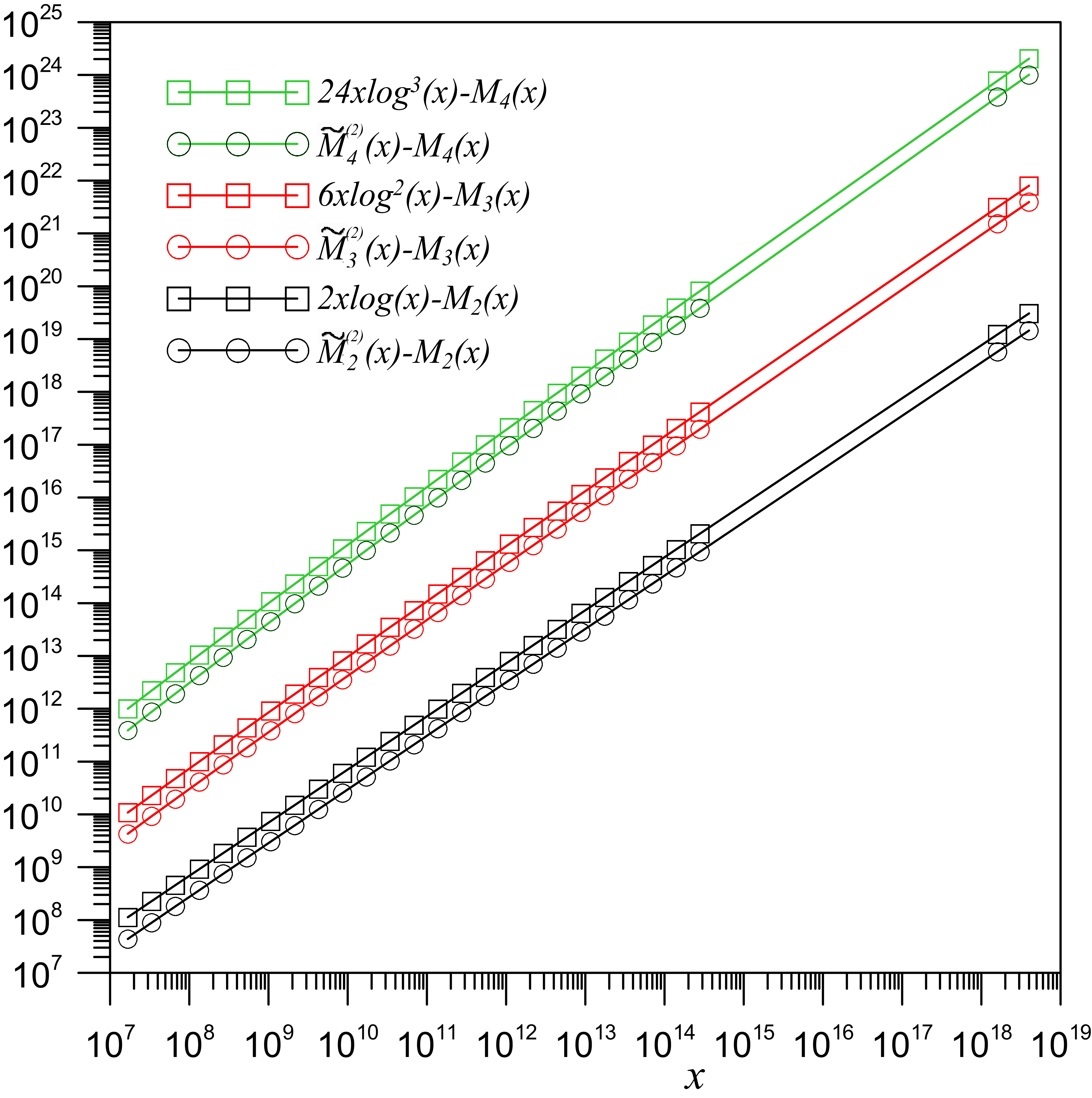} \\  
\caption{\small The plot of  differences between experimental  values of moments $_k(x)$ and calculated from
$\widetilde{M}_k^{(1)}$ and $\widetilde{M}_k^{(2)}$ for $ k=2,3, 4$.
}
\label{Fig-errory}
\end{figure*}

We can try to determine the form of error terms in the formulas  \eqref{prediction1}, \eqref{eq-moment-2}, \eqref{eq-moment-3}
and \eqref{prediction2}.
In   figure  \ref{Fig-errory},  we present plots of the differences of  experimental  values of
moments $M_k(x)$  calculated from the real computer data and appropriate formulas  for  $\widetilde{M}_k^{(i)}(x)$.
All  these plots suggest that  the error term is given by $A_k x^\alpha$, where $\alpha$  is very close to 1 and  the prefactors
$A_k$  increases rapidly with the  order $k$ of  moments.  Because all approximate expressions $\widetilde{M}_k^{(i)}(x)$  give  
values larger  than   experimental  values of moments we write:
\bee
M_k(x)=\widetilde{M}_k^{(i)}(x)-A_k^{(i)}x.
\eee
In the Table IV  we present a sample of coefficients $A_k^{(i)}$ calculated  from the above equation for $x=4\times 10^{18}$,
as then the exponent in power of $x$ is closest to 1.  Thus,   generalizing to non--integer $k$, we  formulate the
{\bf Conjecture}:
\bee
M_k(x) = \frac{ \Gamma(k+1)  x^k}{\pi^{k-1}(x)} + \mathcal{O}(x).          
\label{prediction3}
\eee

\vskip 0.4cm
\begin{center}
{\sf TABLE {\bf IV}}   ~Prefactors  $A_k^{(i)}$  calculated for $x=4\times 10^{18}$.
\begin{tabular}{|c|c|c|c|} \hline
       &  $  A_k^{(1)}  $   &   $ A_k^{(2)}  $ & $ A_k^{(3)}  $   \\ \hline
  $  k=2  $  &      7.674  &       3.624  &       5.624   \\ \hline
  $  k=3 $  &   2003.517  &     985.198  &    1482.890   \\ \hline
  $  k=4  $  &  508697.096  &   252978.305  &   376492.431   \\ \hline
\end{tabular} \\
\end{center}
\vskip 0.4cm

In  paper \cite{Oliveira-2013} on p. 2057 the  authors consider corrections  to \eqref{eq-Heath-Brown}    
given by the series in powers of $1/\log(x)$:
\bee
M_k(x)\approx  k! \log^{k-1}(x) \sum_{n=0}^N \frac{d_{kn}}{\log^n(x)}.
\eee
In this paper the table of values of  $d_{kn}$ obtained from the  least -- square fitting to data  for  $10^{10}<x<4\times 10^{18}$
is given for $N=2$  for $k=2,3, 4$. We have checked that increasing the order $N$  completely changes the values of coefficients
$d_{kn}$, except $d_{k0}\approx 1$,   thus they depend on the order $N$. To  explain   this we notice that the fitting
was done in the very short interval  $(1/\log(4\times 10^{18},   1/\log(10^{10}) = (0.0233\ldots,  0.031\ldots)$.
On such a narrow  interval  each smooth  function  by the Taylor expansion  is a linear function in the first approximation  plus a part of parabola plus a cubic term etc.
In the Taylor expansion of $f(x)$  around  point  $x=a$  coefficients are $f^(n)(a)/n!$   and they does not  change with increasing    
the number of terms.  However  in  \cite{Oliveira-2013}  $d_{kn}$  were determined from the  least --square method.

The correct expansion in powers of  $1/\log(x)$ of formulas  for moments we obtain using the asymptotic series for the logarithmic integral in
\eqref{PNT}   and putting it into ours expressions for moments involving the prime counting function $\pi(x)$.  In this manner we
obtain from \eqref{eq-moment-2} for second  moment:
\bee
\widetilde{M}_2^{(2)}=2x \log(x)\left(1-\frac{2}{\log(x)}-\frac{2}{\log^2(x)}-\frac{3}{\log^3(x)}+\frac{17}{\log^4(x)}+\ldots\right)
\eee
and for third moment:
\bee
 \widetilde{M}_3^{(2)}=3x \log^2(x)\left(1-\frac{4}{\log(x)}+\frac{5}{3\log^2(x)}-\frac{2}{\log^3(x)}+\frac{47}{\log^4(x)}+\ldots \right).
\eee
In general from our conjecture \eqref{prediction2} we get
\bee
\widetilde{M}_k^{(3)}=k!x \log^{k-1}(x)\left(1+\frac{1-k}{\log(x)}+\frac{4-5k+k^2}{2\log^2(x)}+\frac{6-\frac{47 k}{6}+2 k^2-\frac{k^3}{6}}{\log^3(x)} + \ldots\right).
\eee
For the coefficients  $d_{k0}$ in \cite{Oliveira-2013} the values  very close to 1 were  obtained and indeed from above expansions
we have that they are  always 1.

\end{document}